\documentclass{amsart}

\usepackage{amsmath,amsthm,amsfonts,amssymb,bm,graphicx, mathrsfs}

\usepackage
{hyperref}
\hypersetup{colorlinks=true,citecolor=blue,linkcolor=blue,urlcolor=blue,
pdfstartview=FitH }

\theoremstyle{plain}
\newtheorem{theorem}{Theorem}

\newtheorem{lemma}{Lemma}

\newtheorem{cor}[theorem]{Corollary}

\numberwithin{equation}{section}

\theoremstyle{definition}

\newtheorem{remark}{Remark}

\renewcommand{\geq}{\geqslant}
\renewcommand{\leq}{\leqslant}

\newcommand{\changed}[1]{{\color{black} #1}}

\usepackage{calc}
\newsavebox\CBox
\newcommand\hcancel[2][0.5pt]{%
  \changed{\ifmmode\sbox\CBox{$#2$}\else\sbox\CBox{#2}\fi%
  \makebox[0pt][l]{\usebox\CBox}%
  \rule[0.5\ht\CBox-#1/2]{\wd\CBox}{#1}}}

\makeatletter
\DeclareRobustCommand\widecheck[1]{{\mathpalette\@widecheck{#1}}}
\def\@widecheck#1#2{%
    \setbox\z@\hbox{\m@th$#1#2$}%
    \setbox\tw@\hbox{\m@th$#1%
       \widehat{%
          \vrule\@width\z@\@height\ht\z@
          \vrule\@height\z@\@width\wd\z@}$}%
    \dp\tw@-\ht\z@
    \@tempdima\ht\z@ \advance\@tempdima2\ht\tw@ \divide\@tempdima\thr@@
    \setbox\tw@\hbox{%
       \raise\@tempdima\hbox{\scalebox{1}[-1]{\lower\@tempdima\box
\tw@}}}%
    {\ooalign{\box\tw@ \cr \box\z@}}}
\makeatother

\begin{document}

\author{Valentin Blomer}
  
\address{Mathematisches Institut, Bunsenstr. 3-5, 37073 G\"ottingen, Germany} \email{vblomer@math.uni-goettingen.de}

 \title{Spectral summation formula for ${\rm GSp}(4)$ and moments of spinor $L$-functions}

\thanks{First author   supported in part  by the  Volkswagen Foundation and  NSF grant 1128155 while enjoying the hospitality of 
the Institute for Advanced Study. The United States Government is authorized to reproduce and distribute reprints notwithstanding any copyright notation herein.}


\begin{abstract}  We compute the first and second moment of the spinor $L$-function at the central point of Siegel modular forms of large weight $k$ with power saving error term and give applications to non-vanishing. 
  \end{abstract}


\setcounter{tocdepth}{2}  \maketitle 

\maketitle

\section{Introduction}

Spectral summation formulae  like the Selberg trace formula or the Petersson-Bruggeman-Kuznetsov formula belong to the strongest tools in analytic number theory. Their use, however, is almost exclusively restricted to the group ${\rm GL}(2)$, and it is very desirable that they become standard technology also in higher rank situations. 

There are roughly three levels of complexity how spectral summation formulae for a given group $G$ or its associated symmetric space can be applied. The simplest application is Weyl's law (which in sufficient generality may still present formidable difficulties), where no Hecke operators are involved. In fact, this application was Selberg's main motivation to develop the  trace formula in the 1950s.   

The next generation features problems in which some uniformity with respect to Hecke operators is required, but any polynomial bounds suffice and all ``off-diagonal'' terms are estimated trivially. These include in particular the equidistribution of  archimedean and non-archimedean spectral parameters, and the distribution of low-lying zeros of $L$-functions with suitably restricted test functions. Only recently, Matz and Templier \cite{MT}, in a remarkable paper,  succeeded  in establishing    this  in full generality for ${\rm GL}(n)/\Bbb{Q}$. 

The most difficult -- and from the point of  view of analytic number theory the most interesting~-- type of application is when the spectral summation formula develops its full force in the sense that the off-diagonal terms are treated non-trivially and further cancellation is detected using   specific structural features of the formula, for instance on average over Hecke operators. A typical situation is the computation of a moment of $L$-functions in a certain spectrally given family, often with applications to subconvexity and/or non-vanishing. This procedure, which  requires extremely detailed information on the arithmetic/geometric side of the trace formula, is fairly standard for ${\rm GL}(2)$   (in particular for the Kuznetsov formula that may even be applied ``backwards''), but to the author's knowledge a finer trace formula analysis in higher rank 
has only been achieved very recently with the ${\rm GL}(3)$ Kuznetsov formula in \cite{BBM} and \cite{BB}. \\

It is this last type of application that we are concerned with in this paper, in the rank 2 situation of  holomorphic Siegel modular forms. The ``type II'' problems for ${\rm GSp}(4)$ mentioned above have been solved recently  by Kowalski, Saha and Tsimerman \cite{KST}; here we pursue the aim to initiate a detailed analysis of   the corresponding rank 2 Petersson formula to obtain information on spectral averages of spinor $L$-functions for ${\rm GSp}(4)$ in the critical strip with applications to non-vanishing.\footnote{The authors of \cite{KST} allude to the possibility of treating  the off-diagonal terms in the ${\rm GSp}(4)$ Petersson formula non-trivially, but write ``\emph{In our case, the complexity of the analogue expansion for Siegel cusp forms makes this a rather doubtful prospect, at least at the moment.}''}\\
 
We proceed to describe our results. For a general introduction to Siegel modular forms, see \cite{BGHZ, Kl}.  Let   $F  \in S_k^{(2)}$ be a Siegel cusp form of even weight $k$  for the group $\Gamma = {\rm Sp}_4(\Bbb{Z})$ that is an eigenform of the Hecke algebra. This is a function on Siegel's upper half plane $\Bbb{H}_2 = \{Z = X + iY \in \text{Mat}_2(\Bbb{C}) \mid Z = Z^{\top}, Y > 0\}$ equipped with the Petersson inner product
\begin{equation}\label{inner1}
\langle F, G\rangle = \int_{\Gamma \backslash \Bbb{H}_2} F(Z) \overline{G}(Z) (\det Y)^k \frac{dX\, dY}{(\det Y)^{3}}.
\end{equation}
We write the Fourier expansion of $F$ as
$$F(Z) = \sum_{T \in \mathscr{S}} a_F(T) (\det T)^{\frac{k}{2} - \frac{3}{4}} e(\text{tr}(TZ))$$
with (real-valued)  Fourier coefficients $a_F(T)$, where $\mathscr{S}$ is the set of symmetric, positive definite, half-integral matrices $T$ with integral diagonal. We denote by $L(s, F)$ the  spinor $L$-function, normalized so that its critical strip is $0 < \Re s < 1$.  This is a degree 4 $L$-function. With $I$ the 2-by-2 identity matrix, let 
\begin{equation}\label{omega}
\omega_{F, k} :=    \frac{\pi^{1/2} }{4}   (4\pi)^{3-2k} \Gamma(k-3/2)\Gamma(k-2) \frac{ a_F(I) ^2}{\| F \|^2},
\end{equation} and let $\mathcal{B}^{(2)}_k$ denote a Hecke basis of $\mathcal{S}^{(2)}_k$ (which  has  cardinality  $\asymp k^3$ \cite[p.\ 69 and p.\ 123]{Kl}).  
\begin{theorem}\label{thm1} Let $k \geq 6$ be even. Then 
\begin{equation}\label{theo1}
\sum_{F \in \mathcal{B}_k^{(2)}} \omega_{F, k} 
 L(1/2, F) = 2L(1, \chi_{-4}) \log  \frac{k}{4\pi^2}  + 2L'(1, \chi_{-4}) + O(k^{-1}),
 \end{equation}
where 
$\chi_{-4}$ is the non-trivial character modulo 4. 
\end{theorem}

For comparison, Kowalski, Saha and Tsimerman proved
$$\sum_{F \in \mathcal{B}_k^{(2)}} \omega_{F, k} L(s, F) \sim \zeta(s+1/2) L(s+1/2, \chi_{-4}), \quad k \rightarrow \infty$$
for $s \not= 3/2$ (because of possible poles) \emph{outside} the critical strip and without error term. We can even go a step further and compute also the second moment with a power saving error term, where -- with applications in mind -- we allow additional twisting by Dirichlet characters.  This is the main result of this article.

\begin{theorem}\label{thm2} Let $k \geq 6$ be even, $\varepsilon > 0$. Let $q_1, q_2$ be two coprime\footnote{The coprimality assumption is for technical convenience only; equation \eqref{theo2} remains true for any real primitive characters.} fundamental discriminants (possibly 1).  Then 
\begin{equation}\label{theo2}
\sum_{F \in \mathcal{B}_k^{(2)}} \omega_{F, k} 
L(1/2, F \times \chi_{q_1})L(1/2, F \times \chi_{q_2}) = \text{{\rm  main term }}+ O_{q_1, q_2}(k^{-1/2 + \varepsilon}),
\end{equation}
where the main term is the residue at $s=t=0$ of the expression \eqref{integrand}. 

In particular, if $q_1 = q_2 = 1$, the main term equals
\begin{equation}\label{a1}
\frac{4}{3} L(1, \chi_{-4})^2 P_3\left(\log k\right)
\end{equation}
for a certain monic polynomial $P_3$ of degree $3$ depending on $q_1, q_2$. 

If $\{q_1, q_2\} \in \{1, -4\}$,   the main term equals 
\begin{equation}\label{a2}
   2L(1, \chi_{-4})^3 P_2\left(\log k\right) 
   \end{equation}
for a certain monic polynomial $P_2$ of degree $2$ depending on $q_1, q_2$. 

If $q_1, q_2$ are two coprime integers different from $1$ and $-4$, the main term equals
\begin{equation}\label{a3}
4 L(1, \chi_{q_1})L(1, \chi_{-4q_1}) L(1, \chi_{q_2})L(1, \chi_{-4q_2})L(1, \chi_{q_1q_2}). 
\end{equation}
The error term in \eqref{theo2} depends polynomially on $q_1, q_2$. 
\end{theorem}

If $q_1 = q_2 = 1$, the main term contains  an off-diagonal contribution coming from Klooster\-man sums and Bessel functions in the Petersson formula. It requires a very careful  analysis    to extract the relevant portion of the main term, which is based both on subtle properties of special functions and an interesting analysis of  symplectic exponential sums (see Section \ref{sec5}). It is a very pleasing structural fact that this extra off-diagonal main term matches precisely a certain polar contribution of the diagonal term. This requires some non-trivial manipulation and is a good way to double-check the somewhat intricate computation\footnote{It turns out that the Fourier expansion of Siegel Poincar\'e series in \cite{Ki} contains a numerical inaccuracy that was found in this way, cf.\ Remark \ref{rem} below.}. Roughly speaking, the off-diagonal term in the Petersson formula is a sum over integral matrices. Special features of (integrals of products of) Bessel functions, exploited in Lemma \ref{psi}, imply that after Poisson summation only matrices in ${\rm GO}_2(\Bbb{Z}) =  \Bbb{R}_{>0} \cdot  O(2) \cap \text{Mat}_2(\Bbb{Z})$ survive, and modulo automorphisms this can be identified as a semigroup with the non-zero integral ideals of $\Bbb{Q}(i)$. This brings us to Dedekind zeta-functions, and after applying their functional equation we recognize an earlier polar term. \\

The ``harmonic'' weights $\omega_{F, k}$ are natural from the point of view of spectral summation formulae of Petersson's type. They are of size $k^{-3}$ on average since
\begin{equation}\label{av}
  \sum_{F \in \mathcal{B}_k^{(2)}} \omega_{F, k}  = 1 + O(e^{-k})
\end{equation}
(cf.\ also \cite[Theorem 4.11]{DPSS} for $N=1$, but notice the much stronger error term in the $k$-aspect in \eqref{av}), but they carry very different, and probably much more complicated, arithmetic information than the corresponding harmonic weights
$$(4\pi)^{1-k} \Gamma(k-1) \frac{a_f(1)^2}{\| f \|^2}$$
  for elliptic modular forms $f \in S_{k}$. While the latter are (only) related to  $L$-values $1/L(1, \text{sym}^2 f)$ at the edge of the critical strip, B\"ocherer \cite{Bo} made a remarkable conjecture that $\omega_{F, k}$ should (in addition) be related to \emph{central} $L$-values.   This can be seen as generalization of Waldspurger's theorem. We refer to \cite[Chapter 1, in particular Conjecture 1.10]{FS} for some enlightening discussion. In particular, it is not even known if $\omega_{F, k}$ can be zero, and if so, how often this can happen.  Recently, a very precise version of B\"ocherer's conjecture was put forward in the beautiful paper \cite[Conjecture 1.2]{DPSS}, which states that 
for non-Saito-Kurokawa lifts 
\begin{equation}\label{conj}
\omega_{F, k}  = \frac{64 \pi^6 \Gamma(2k-4)}{\Gamma(2k-1)} \frac{L(1/2, F) L(1/2, F \times \chi_{-4})}{L(1, \pi_F, \text{Ad})},   
\end{equation}
where $\pi_F$ is the automorphic representation associated with $F$ and $L(s, \pi_F, \text{Ad})$ is a degree 10 $L$-function.\footnote{By \cite[Theorem 5.2.1]{PSS} it is known that  $L(1, \pi_F, \text{Ad}) \not= 0$, but it seems that no good lower bounds are available. Upper bounds follow from \cite{Li}.}
 If this is true, then Theorem \ref{thm2} really evaluates a \emph{fourth moment} of central values, and hence a degree 16 $L$-function!  

Theorem \ref{thm2} gives further evidence towards the B\"ocherer conjecture as it shows unconditionally that there is stronger correlation between $\omega_{F, k}$ and $L(1/2, F \times \chi_{q_1}) L(1/2, F \times \chi_{q_2})$ if $q_1, q_2\in  \{1, -4\}$.

 If   
$F \in S_k^{(2)}$ is a Saito-Kurokawa lift coming   from an elliptic modular form $f \in S_{2k-2}$, then a variant of \eqref{conj} is a theorem, and we have (see e.g. \cite[Section 4.4]{DPSS}) 
$$\omega_{F, k} \asymp \frac{1}{k^3} \frac{ L(1/2, f \times \chi_{-4})}{L(1, \text{sym}^2 f)}.$$ 
Combining this with \eqref{saito} below, we see that the contribution of lifts to the left hand side of \eqref{theo1} and \eqref{theo2} is very small, in fact $O(k^{-1/2+\varepsilon})$ using only the convexity bound for the central $L$-values, but much better bounds could be obtained. 

By the Cauchy-Schwarz inequality we obtain  from Theorem \ref{thm1} and \eqref{a1} immediately

\begin{cor}\label{cor} Let $k \geq 6$ be sufficiently large and let $(\mathcal{B}_k^{(2)})^{\ast}$ denote a   Hecke basis of the space orthogonal to Saito-Kurokawa lifts. Then  
$$ \sum_{\substack{F \in( \mathcal{B}_k^{(2)})^{\ast} \\ L(1/2, F) \not= 0} } \omega_{F, k}   \gg (\log k)^{-1}.$$
In particular, there exist (generic) cusp forms $F \in (\mathcal{B}_k^{(2)})^{\ast}$ with   $a_F(I) L(1/2, F) \not= 0$. 
\end{cor}

Somewhat similar in spirit is the non-vanishing result of \cite{DK}  for the Koecher-Maa{\ss} $L$-function associated to $F$ (having a functional equation, but no Euler product) in the critical strip, but outside the critical line. \\

From \eqref{a2} and \eqref{av} and the Cauchy-Schwarz inequality we obtain

\begin{cor}\label{cor1} Assume that \eqref{conj} holds. Then
$$ \sum_{\substack{F \in (\mathcal{B}_k^{(2)})^{\ast} \\ \omega_{F, k} \not= 0} }  \frac{1}{L(1, \pi_F, \text{{\rm Ad}})} \gg \frac{k^3}{(\log k)^{2}}.$$
In particular, if $L(s, \pi_F, \text{{\rm Ad}})$ has no zeros in $|s-1| \ll k^{-\varepsilon}$, then $k^{3-\varepsilon}$ members $F \in  (\mathcal{B}_k^{(2)})^{\ast}$ satisfy $\omega_{F, k}\not=0$  (and hence $  L(1/2, F) L(1/2, F \times \chi_{-4}) \not= 0$). 
\end{cor}

Finally, from \eqref{a3} we obtain the following \emph{quadruple} non-vanishing result. 
\begin{cor}\label{cor2} Assume that \eqref{conj} holds. Let   $q_1, q_2$ be any two coprime fundamental discriminants and let $k$ be sufficiently large. Then there exists $F \in  (\mathcal{B}_k^{(2)})^{\ast}$ such that
$$L(1/2, F) L(1/2, F \times \chi_{-4}) L(1/2, F \times \chi_{q_1})L(1/2, F \times \chi_{q_2}) \not= 0.$$\\
\end{cor}

The above discussion shows that central values of $L$-functions of cohomological type for ${\rm GSp}(4)$ belong to the most fascinating arithmetic-analytic objects. 
It is therefore of interest to investigate  their analytic properties, and the above results seem to provide the first analytic properties of symplectic $L$-functions inside the critical strip. 
The technology developed in this paper is   capable of several extensions of which we mention four:

(a) It is possible to include more general weights $\omega_{F, k, \Delta, \Lambda}$ for a negative fundamental discriminant $\Delta$ and a class group character $\Lambda \in  \widehat{{\rm Cl}}(\Delta)$, where $|a_F(I)|^2$ in \eqref{omega}  is replaced with 
$$ A_F(\Delta)  := \frac{1}{h(\Delta)} \Bigl|\sum_{J \in {\rm Cl}(\Delta)} \Lambda(J) a_F(J)\Bigr|^2$$
(where we identify an element $J$ in the class group with the matrix of an associated integral quadratic form). If the class number of $\Bbb{Q}(\sqrt{\Delta})$ is 1, this requires only notational changes.

 (b) Because of the power saving error terms    in Theorems \ref{thm1} and \ref{thm2}, one can insert   additional Dirichlet polynomials such as amplifiers or mollifiers. 
 
 (c) One can average over $k$ in some dyadic interval $K \leq k \leq 2K$, which may enable one to treat higher moments (e.g.\ a fourth moment might be within reach, although it is certainly a challenge). 
 
 (d) It is also  possible to treat the level aspect. This is maybe the most interesting variation, as it gives rise to richer families of $L$-functions associated to algebraic modular forms (which also contain Yoshida lifts of certain pairs of elliptic cusp forms of weight 2 and $2k-2$). There are two natural candidates for congruence subgroups, the paramodular subgroup and the Siegel subgroup. We refer to  \cite{Sch} for progress towards a newform theory in these cases. As a relatively simple sample result in this direction we have the following.
 \begin{theorem}\label{thm3}  Let $N \equiv 3 \, (\text{{\rm mod }} 4)$ be a large prime, and let $\Gamma = \Gamma^{(2)}_0(N) \subseteq {\rm Sp}_4(\Bbb{Z})$ be the Siegel congruence subgroup consisting of matrices $\left(\begin{smallmatrix} A & B\\ C & D\end{smallmatrix}\right)$ with $C \equiv 0$ {\rm (mod }$N)$. Let $k \geq 6 $ be even, and let  $\mathcal{B}_k^{(2)}(N)$ denote  a corresponding Hecke basis. Then\footnote{The error term can be improved considerably if desired.}
 $$\sum_{F \in \mathcal{B}_k^{(2)}} \omega_{F, k}   L(1/2, F) = 2L(1, \chi_{-4}) \log N + C(k) +     O_k(N^{-1})$$
 for an explicitly given constant $C(k)$. In particular, for   sufficiently large $N$ not all central values can vanish. 
 \end{theorem}


\textbf{Notation and conventions.} We use the usual $\varepsilon$-convention, and all implied constants may depend on $\varepsilon$. We refer to a quantity as \emph{negligible} if it is $\ll k^{-100}$. For notational simplicity we write $\ell := k - 3/2$. We write $[., .]$ for the positive least common multiple of two non-zero integers.  \\

\textbf{Acknowledgement.} The author is very grateful to Abhishek Saha for useful comments and suggestions.

\section{The spinor $L$-function}

For a Siegel cusp form $F  \in S_k^{(2)}$ of even weight $k$ that is an eigenform of the Hecke algebra with local parameters $\alpha_{0, p}$, $\alpha_{1, p}$, $\alpha_{2, p}$ (satisfying $\alpha_{0, p}^2 \alpha_{1, p} \alpha_{2,p} = 1$) at primes $p$, the spinor $L$-function is defined by a degree 4  Euler product
$$L(s, F) = \prod_{p} \left( 1 - \frac{\alpha_{0, p}}{p^s}\right)^{-1} \left( 1 - \frac{\alpha_{0, p}\alpha_{1, p}}{p^s}\right)^{-1} \left( 1 - \frac{\alpha_{0, p}\alpha_{2, p}}{p^s}\right)^{-1} \left( 1 - \frac{\alpha_{0, p}\alpha_{1, p}\alpha_{2, p}}{p^s}\right)^{-1} $$
for $\Re s$ sufficiently large. Its meromorphic continuation and functional equation was proved by Andrianov \cite[Theorem 3.1.1]{An}:  
\begin{equation}\label{fe}
\Lambda(s, F) = L_{\infty}(s, F) L(s, F) =  \Gamma_{\Bbb{C}}(s+1/2)\Gamma_{\Bbb{C}}(s + k - 3/2) L(s, F) = \Lambda(1 - s, F),
\end{equation}
where as usual $\Gamma_{\Bbb{C}}(s) = 2(2\pi)^{-s} \Gamma(s)$. 
Notice that in contrast to most of the classical literature on Siegel modular forms we normalize all $L$-functions to have $0 < \Re s < 1$ as the critical strip. This is very convenient and corresponds to a linear shift $s \mapsto s + k - 3/2$ in comparison with \cite{An} and many other sources. We also normalize the  Hecke eigenvalues $\lambda_F(m)$ of $F$ accordingly. The Dirichlet series expansion of $L(s, F)$ is given by \cite[p.\ 69]{An}
$$L(s, F) = \zeta(2s+1) \sum_m \frac{\lambda_F(m)}{m^s}.$$

The space $S_k^{(2)}$ contains a subspace of (Saito-Kurokawa) lifts from elliptic Hecke cusp forms $f \in S_{2k-2}$ of weight $2k-2$. For such lifts $F$ corresponding to $f$, we have 
\begin{equation}\label{saito}
L(s, F) = \zeta(s + 1/2) \zeta(s-1/2) L(s, f),
\end{equation}
in particular these functions have a pole at $s = 3/2$ (but no pole at $s = 1/2$ since $2k-2 \equiv 2$ (mod 4), so $L(1/2, f) = 0$ for root number reasons). Except for the pole at $s=3/2$ for Saito-Kurokawa lifts, $L(s, F)$ is entire. If $F$ is not a lift, then the Ramanujan conjecture holds by a deep result of Weissauer \cite{We} (that we will not use in this paper), which in our normalization states $\lambda_F(m) \ll m^{\varepsilon}$. 

Let $q$ be a fundamental discriminant (possibly 1). The twisted $L$-functions $L(s, F\times \chi_q)$ is given by twisting its Dirichlet series expansion with the character $\chi_q$. It satisfies  a similar self-dual functional equation 
\begin{equation}\label{fe1}
\Lambda(s, F \times \chi_{q})  = |q|^{2s}  \Gamma_{\Bbb{C}}(s+1/2)\Gamma_{\Bbb{C}}(s + k - 3/2) L(s, F \times \chi_{q}) = \Lambda(1 - s, F \times \chi_{q}),
\end{equation} 
which was recently established in  full generality by Krieg and Raum \cite{KR}. 

As usual in higher rank, the Fourier coefficients cannot easily be recovered from the Hecke eigenvalues, but conversely the Hecke eigenvalues can be written rather simply in terms of Fourier coefficients. 
An important ingredient in  Andrianov's proof of the functional equation is the explicit formula \cite[Theorem 2.4.1]{An}, a special case of which is
\begin{equation}\label{andr}
L(s+1/2, \chi_q)L(s+1/2, \chi_{-4q} )  \sum_m \frac{a_F(mI) \chi_{q}(m)}{m^s} = L(s, F \times \chi_{q})  a_F(I), 
\end{equation}
see \cite[Theorem 4.3.16]{An2} with $l = a = 1$, $\eta = \chi = \text{triv}$. We denote by
$$r(n) = r_{q}(n) = \frac{\chi_q(n)}{n^{1/2}} \sum_{d \mid n}  \chi_{-4}(d)$$
the Dirichlet coefficients of $L(s+1/2, \chi_q)L(s+1/2, \chi_{-4q})$. In particular, if $q=1$, the latter is the Dedekind zeta function $\zeta_{\Bbb{Q}(i)}(s+1/2)$.

From \eqref{fe} and \eqref{andr} we obtain in a standard fashion \cite[Theorem 5.3]{IK} an approximate functional equation 
\begin{equation}\label{approx}
\begin{split}
&a_F(I) L(1/2, F \times \chi_q) = 2 \sum_{n, m} \frac{r(n) a_F(mI) \chi_q(m)}{(nm)^{1/2}} W(nm/|q|^2),\\& W(x) = \frac{1}{2\pi i} \int_{(2)} \frac{L_{\infty}(s + 1/2)}{L_{\infty}(1/2)} (1-s^2)x^{-s} \frac{ds}{s}.
\end{split}
\end{equation}
The factor $(1-s^2)$ was inserted to make sure that the integrand vanishes at $s=1$ to counteract the pole of Saito-Kurokawa  lifts. Any even, polynomially bounded, holomorphic  function $G(s)$ satisfying $G(1) = 0$, $G(0) = 1$ would serve equally well, and we will see that the factor $1-s^2$ will be irrelevant in all forthcoming residue computations.  The integral is rapidly converging; differentiating under the integral sign and  using Stirling's formula  it is easy to see 
by shifting the contour to either $\Re s = -1/2$ (say) or to $\Re s = A$ that 
\begin{equation}\label{boundW}
x^j W^{(j)}(x) \ll_{j, A} (1 + x/k)^{-A}
\end{equation}
 for any $A > 0$ and any $j \in \Bbb{N}_0$.

\section{The Petersson formula}\label{sec3}

Our main tool in this paper is a spectral summation formula of ``Petersson type'', which in the context of Siegel modular forms can be proved in the same way as in the classical case, namely by computing the inner product of two Poincar\'e series. The relevant Fourier expansion of the Poincar\'e series has been worked out by Kitaoka \cite{Ki}. We quote his results and introduce some notation.
For $Q, T \in \mathscr{S}$ and an invertible matrix $C \in \text{Mat}_2(\Bbb{Z})$ denote by 
\begin{equation}\label{klooster}
K(Q, T; C) = \sum e\left(\text{tr}(A C^{-1} Q + C^{-1} D T)\right)
\end{equation}
the ``Kloosterman sum'', where the sum is taken over matrices $\left(\begin{smallmatrix} A & \ast \\ C & D\end{smallmatrix}\right) \in {\rm Sp}_4( \Bbb{Z})$ for a given value of $C$ in a system $X(C)$ of representatives for $\Gamma_{\infty} \backslash {\rm Sp}_4( \Bbb{Z}) /\Gamma_{\infty}$ with $\Gamma_{\infty} = \left\{ \left(\begin{smallmatrix} I & X \\  & I\end{smallmatrix}\right)  \mid X = X^{\top}\right\}.$ The cardinality of $X(C)$ depends only on the elementary divisors of $C$ (since $K(Q, T; U^{-1}CV^{-1}) = K(UQ U^{\top}, V^{\top} T V; C)$ for $U, V \in {\rm GL}_2(\Bbb{Z})$), and for $C = \left(\begin{smallmatrix} c_1 & \\ & c_1c_2\end{smallmatrix}\right)$ one sees easily that $|X(C)| \leq c_1^3 c_2 \leq |\det C|^{3/2}$, which is a trivial upper bound for $K(Q, T; C)$.

For a real, diagonalizable matrix $P$ with positive eigenvalues $s_1^2, s_2^2$ ($s_1, s_2 > 0$) we write $$\mathcal{J}_{\ell}(P) := \int_0^{\pi/2} J_{\ell}(4\pi s_1 \sin \theta)J_{\ell}(4\pi s_2 \sin \theta)\sin\theta\, d\theta.$$
For two matrices $P = \left(\begin{smallmatrix} p_1& p_2/2\\ p_2/2 & p_4\end{smallmatrix}\right) \in \mathscr{S}$, $S = \left(\begin{smallmatrix} s_1& s_2/2\\ s_2/2 & s_4\end{smallmatrix}\right) \in \mathscr{S}$ and $c \in \Bbb{N}$ we define another Kloosterman sum
$$H^{\pm}(P, S; c) =  \delta_{s_4 = p_4} \underset{d_1 \, \text{mod }c}{\left.\sum\right.^{\ast}} \sum_{d_2 \, \text{mod }c} e \left( \frac{\overline{d_1} s_4 d_2^2 \mp \overline{d_1} p_2d_2 + s_2d_2 + \overline{d_1}p_1 + d_1 s_1 }{c} \mp \frac{p_2s_2}{2 cs_4}\right).$$
We note in passing that this sum -- essentially a Sali\'e sum -- comes up also in the Fourier expansion of  Jacobi Poincar\'e series of index $s_4 = p_4$ \cite[p.\ 519]{GKZ}. We have the trivial bound $|H^{\pm}(P, S; c)| \leq c^2.$ 

To derive the Petersson formula, we define for $Q \in \mathscr{S}$ a Poincar\'e series
$$P_Q(Z) = \sum_{\gamma \in \Gamma_{\infty} \backslash {\rm Sp}_4(\Bbb{Z})} \det (J(\gamma, Z))^{-k} e(\text{tr}(Q \gamma Z)) = \sum_{T \in \mathscr{S}} h_Q(T) (\det T)^{\frac{k}{2} - \frac{3}{4}}  e(\text{tr}(TZ)),$$
say, where $J(\gamma, Z) = CZ + D$ for $\gamma = \left(\begin{smallmatrix} A &B \\ C & D\end{smallmatrix}\right)$.   Then $$\langle F, P_Q \rangle = 8 c_k (\det Q)^{-\frac{k}{2} + \frac{3}{4}}a_F(Q), \quad c_k =  \frac{1}{4}   \pi^{1/2}   (4\pi)^{3-2k} \Gamma(k-3/2)\Gamma(k-2)$$ (see \cite[(3.1.1)]{KST}, but note our different normalization of the Fourier coefficients)  
for $F \in S_k^{(2)}$. 
Computing $\langle P_T, P_Q\rangle$ for $T, Q \in \mathscr{S}$, we conclude
\begin{equation}\label{inner}
8 c_k \left(\frac{\det T}{\det Q}\right)^{\frac{k}{2} - \frac{3}{4}} \sum_{F \in \mathcal{B}_k^{(2)}} \frac{a_F(T) a_F(Q)}{\| F \|^2} = h_Q(T) (\det T)^{\frac{k}{2} - \frac{3}{4}}.
\end{equation}
The Fourier coefficients $ h_Q(T) (\det T)^{\frac{k}{2} - \frac{3}{4}}$ have been computed by Kitaoka \cite[Sections 2-4]{Ki} (see also \cite[pp.\ 356-358]{KST}), and we quote the following for convenient reference.
\begin{lemma}\label{1}
For $T, Q\in \mathscr{S}$ and even $k \geq 6$ we have 
\begin{displaymath}
\begin{split}
& h_Q(T)   (\det T)^{\frac{k}{2} - \frac{3}{4}}  = \delta_{Q \sim T} \#\text{{\rm Aut}}(T)\\
& +  \left(\frac{\det T}{\det Q}\right)^{\frac{k}{2} - \frac{3}{4}} \sum_{\pm} \sum_{s, c\geq 1} \sum_{U, V} \frac{(-1)^{k/2} \sqrt{2} \pi}{c^{3/2}s^{1/2}} H^{\pm} (UQU^{\top}, V^{-1} T V^{-\top} , c)J_{\ell}  \left(\frac{4\pi \sqrt{\det(TQ)}}{c s}\right)\\
& + 8\pi^2 \left(\frac{\det T}{\det Q}\right)^{\frac{k}{2} - \frac{3}{4}}  \sum_{\det C \not= 0} \frac{K(Q, T; C)}{|\det C|^{3/2}}  \mathcal{J}_{\ell}(T C^{-1} Q C^{-\top}),
\end{split}
\end{displaymath}
where the sum over $U, V \in {\rm GL}_2(\Bbb{Z})$ in the second term on the right hand side is over matrices 
$$U  = \left(\begin{matrix} * & *\\ u_3 & u_4 \end{matrix}\right)/\{\pm 1\},  \quad  
V  = \left(\begin{matrix} v_1 & *\\ v_3 &* \end{matrix}\right), \quad  (u_3\, u_4) Q \left(\begin{matrix} u_3\\ u_4\end{matrix}\right) = (-v_3\, v_1) T \left(\begin{matrix} -v_3\\ v_1\end{matrix}\right) = s,
$$ 
$Q \sim T$ means equivalence in the sense of quadratic forms and $\text{{\rm Aut}}(T) = \{U \in {\rm GL}_2(\Bbb{Z}) \mid U^{\top} T U = T\}$. The sums are absolutely convergent for $k \geq 6$. 
\end{lemma}

\begin{remark}\label{rem} Kitaoka \cite[p.\ 166]{Ki} has the constant $1/2\pi^4$ instead of $8\pi^2$ in the last line, but this turns out to be incorrect: the factor $2(2\pi)^{-3}$ in the third last display on p.\ 165 belongs to the other side of the equation, as can be seen by comparing with p.\ 478 for the choice of measure and p.\ 486 for the definition of $A_{\delta}(M)$ in \cite{He}. (The second last display in \cite[p.\ 165]{Ki}   coincides with \cite[p.\ 517]{He} and is correct.) 
While our results do not depend on the value of the constants, their values are responsible for the matching of two terms in Section \ref{sec6}, which is a structurally important feature. 
\end{remark}

Following \cite{Ki} and \cite{KST}, we refer to the second term on the right hand side as the rank 1 case and to the third term as the rank 2 case. Notice that $T C^{-1} Q C^{-\top}$ is a product of positive symmetric matrices and hence diagonalizable (not necessarily symmetric) with positive eigenvalues, so that $\mathcal{J}_{\ell}(T C^{-1} Q C^{-\top})$ makes sense.  We have the following simple lemma.
\begin{lemma}\label{2}  
For positive definite matrices $T, Q$ with largest eigenvalues $\lambda_T$, $\lambda_Q$, the smallest eigenvalue of $T C^{-1} Q C^{-\top}$ is   $\ll \lambda_T \lambda_Q  \| C \|^{-2}$.
\end{lemma} 

\noindent \textbf{Proof.} The smallest eigenvalue $\lambda_{\min}$ of $T C^{-1} Q C^{-\top}$ is   the inverse of the largest eigenvalue of $(T C^{-1} Q C^{-\top})^{-1}$, which is conjugated to $D^{\top} D$ with $D = Q^{-1/2} C T^{-1/2}$. Hence
$$\lambda_{\min} \ll \|D \|^{-2} = \| Q^{-1/2} C T^{-1/2} \|^{-2}  = \frac{\| Q^{1/2} \|^2 \| T^{1/2} \|^2 }{(\| Q^{1/2} \| \|  Q^{-1/2} C T^{-1/2}  \| \| T ^{1/2} \|)^2} \leq \frac{\| Q^{1/2} \|^2 \| T^{1/2} \|^2}{ \| C \|^2}.$$ 
\\

To check absolute convergence in Lemma \ref{1}  for $k \geq 6$ (for $k \leq 8$ the space $S_k^{(2)}$ is $\{0\}$, so this is no loss of generality), we notice that  the number of representations of $s$ by an integral  positive definite quadratic form is $O(s^{\varepsilon})$ and $J_k(x) \ll_k  x^k$, hence the rank 1 term is 
$$\ll_{T, Q, k} \sum_{c, s} \frac{c^{1/2}}{s^{1/2 - \varepsilon}} \frac{1}{(cs)^{k - 3/2}} < \infty,$$
while the rank 2 term is 
$$\ll_{T, Q, k} \sum_{\det C \not= 0} \frac{1}{\| C \|^{k - 3/2}} \ll \sum_{c \geq 1} \frac{1}{c^{k - 3 - 3/2}} < \infty.$$ 

\begin{remark}\label{rem2} For the (Siegel type) congruence subgroup $\Gamma^{(2)}_0(N)\subseteq {\rm Sp}_4(\Bbb{Z})$, the only modifications in Lemma \ref{1} are the additional congruence conditions $N \mid c$ in the rank 1 term and $N \mid C$ in the rank 2 term, see \cite{CKM}.  
\end{remark}

\section{Interlude: Bessel functions}

In this section we compile all necessary information needed on the $J_k$ function and integrals thereof, where we think of $k$ as being large.  We start with the two simple uniform bounds
\begin{equation}\label{bessel1}
  J_k(x) \ll 1
\end{equation}
and 
\begin{equation}\label{bessel2}
  J_k(x) \ll \left(\frac{x}{k}\right)^k,  
\end{equation}
valid for $x > 0$, $k > 2$, which follow immediately from the integral representations \cite[8.414.13]{GR} and \cite[8.414.4]{GR}, respectively.  We also have the more refined uniform upper bound
\begin{equation}\label{bessel3}
  J_k(x) \ll \min(k^{-1/3}, |x^2 - k^2|^{-1/4}),
\end{equation}
which follows from Olver's uniform expansion \cite[(4.24)]{Ol}. Coupled with an asymptotic expansion of the Airy function at large negative arguments \cite[(4.07)]{Ol2}, this gives in particular
\begin{equation}\label{unif}
\begin{split}
& J_k(x) =  \frac{\exp(i\psi(x, k) ) F_k^{+}(x) + \exp(-i\psi(x, k) ) F^{-}_k(x)}{(x^2 - k^2)^{1/4}} + O(k^{-200}), \\
& \psi(x, k) =  \sqrt{x^2 - k^2} - k \arctan(\sqrt{(x/k)^2 - 1})
\end{split}
\end{equation}
for $x \geq 2k$ (in fact $x \geq k + k^{1/3+\varepsilon}$ would suffice), where $F_k^{\pm}$ are smooth non-oscillating functions  satisfying the uniform bounds $x^j (F_k^{\pm})^{(j)}(x) \ll_j 1$ for all $j \in \Bbb{N}_0$.  

We   record  the differentiation rule \cite[8.471.2]{GR} 
\begin{equation}\label{diff}
J_k'(x) = \frac{1}{2} \left(J_{k-1}(x) - J_{k+1}(x)\right)
\end{equation}
and the Mellin formula \cite[6.574.2]{GR}
\begin{equation}\label{mellin}
\int_0^{\infty} J_k(1/x)^2 x^{s-1} dx = \frac{ \Gamma(k-s/2)\Gamma((1+s)/2)}{2\pi^{1/2}  \Gamma(1+k+s/2)\Gamma(1 + s/2)}
\end{equation}
for   $2k > \Re s > -1$.  A fundamental role is played by the formula (\cite[2.12.20]{PBM}\footnote{Note that a factor $\pi$ is missing in \cite[2.12.20.5]{PBM}  in comparison with \cite[2.12.20.2/3]{PBM}; the same oversight occurs in \cite[1.6.37]{Ob}. For our purposes, of course, the numerical constant is irrelevant.})
\begin{equation}\label{Iwali}
2\Re \Biggl( e\left(- \frac{k+1}{4}\right) \int_0^{\infty} e\left((\alpha+\beta)z + \frac{\gamma}{z}\right)J_k(4\pi \sqrt{\alpha \beta}z )\frac{dz}{z}\Biggr)    = 2\pi J_k(4\pi \sqrt{\alpha \gamma}) J_k(4\pi \sqrt{\beta \gamma}) 
\end{equation}
for $\alpha, \beta, \gamma > 0$ that was used in a very different context in \cite[(A.9)]{IL}. 

The central aim in this section is to understand the Fourier integral
\begin{equation}\label{defPsi}
\begin{split}
&\Psi(C; h_1, h_2) = \Psi_{n_1, n_2, q_1, q_2}(C; h_1, h_2)\\
& := \int_{\Bbb{R}} \int_{\Bbb{R}} \frac{W(n_1x_1/|q_1|^2)W(n_2x_2/|q_2|^2)}{\sqrt{x_1x_2}} \mathcal{J}_{\ell}\left(x_1x_2 (C^{\top} C)^{-1}\right) e\left( \frac{x_1h_1 + x_2h_2}{|\det C|}\right)dx_1dx_2,
\end{split}
\end{equation}
where $n_1, n_2 \in \Bbb{N}$, $h_1, h_2 \in \Bbb{Z}$, $q_1, q_2 \in \Bbb{Z} \setminus \{0\}$,   $W$ is given by \eqref{approx} and $C\in \text{Mat}_2(\Bbb{Z})$ with non-zero determinant.  We denote by $ s_1^2, s_2^2$ ($s_1, s_2 > 0$) the eigenvalues of $(C^{\top} C)^{-1}$. Clearly there is no oscillation in the integral if $s_1 = s_2$ and $h_1 = h_2 = 0$. Our next lemma shows that in all other cases $\Psi(C; h_1, h_2)$  is small in the ranges of $h_1, h_2, C$ that will be of interest for us. 
\begin{lemma}\label{psi}
Suppose 
$h_1, h_2, \| C \| \ll k^{\varepsilon}.$ Then $$\Psi(C; h_1, h_2) \ll_{q_1, q_2} k^{-1/2 +\varepsilon}$$
(uniformly in $n_1, n_2$) unless $h_1 = h_2 = 0$ and $s_1 = s_2$. 
\end{lemma} 

\textbf{Proof.}  Since $C$ is integral and $\| C \| \ll k^{\varepsilon}$,   we have $ k^{-\varepsilon} \ll s_1, s_2 \ll k^{\varepsilon}$. The integrality of $C$ also implies the eigenvalues of $C^{\top} C$ are either identical or differ at least by 1, hence $s_1, s_2$ are identical or differ by at least $\gg k^{-\varepsilon}$. In the following all implied constants may depend on $q_1, q_2$. 

Now \eqref{bessel2} implies that the integral defining $\Psi(C; h_1, h_2)$ is negligible in the range $x_1, x_2 \ll k^{1-\varepsilon}$, and hence $n_1, n_2 \ll k^{\varepsilon}$ by \eqref{boundW}. We remember this by inserting two smooth, non-negative functions $v$ that are 1 for $x \geq k^{1-\varepsilon}$ and 0 for $x \leq \frac{1}{2} k^{1-\varepsilon}$.  Using \eqref{Iwali}, we can write  $\Psi(C; h_1, h_2)$ as a sum of two terms of the form 
\begin{equation}\label{mult}
\begin{split}
   \int_{\Bbb{R}} \int_{\Bbb{R}} \int_0^{\pi/2}  \int_0^{\infty} & v(x_1) v(x_2) \frac{W(n_1x_1/|q_1|^2)W(n_2x_2/|q_2|^2)}{\sqrt{x_1x_2}}   e\left(\pm \Bigl((s_1^2+s_2^2)z + \frac{x_1x_2 \sin^2\theta}{z}\Bigr)\right)\\
   &\times J_{\ell}\left(4\pi s_1s_2z \right)\frac{dz}{z} \,  \sin\theta \, d\theta\, 
 e\left( \frac{x_1h_1 + x_2h_2}{|\det C|}\right)dx_1dx_ 2 + O(k^{-100}). 
 \end{split}
\end{equation}
This multiple integral is absolutely convergent. We distinguish   two cases depending on the size of $z$.  Let $V$ be a smooth, non-negative function that is 1 for $z \geq k^{1+\delta}$ and 0 for $z \leq \frac{1}{2} k^{1+\delta}$ for some small $\delta > 0$.  

We first insert $V(z)$ and treat the portion where $z$ is large. Integrating by parts sufficiently often with respect to $x_1, x_2$ shows that the integral is negligible unless $h_1 = h_2 = 0$ (since $\det C \ll k^{\varepsilon}$). Moreover, we can insert the uniform asymptotic expansion \eqref{unif} for the Bessel function with a negligible error term. Then we can write the $z$-integral as two terms of the form
\begin{displaymath}
\begin{split}
\int_0^{\infty} \frac{V( z)}{((4\pi  s_1s_2z)^2 - k^2)^{1/4}} F^{\epsilon_1} (4\pi  s_1s_2 z) e\left(\frac{x_1x_2 \sin^2\theta}{z} -\frac{\epsilon_1 k}{2\pi} \arctan \sqrt{(4\pi  s_1s_2z/k)^2 - 1}\right)\\ \times e\left( -\frac{\epsilon_1 k}{2\pi (\sqrt{(4\pi  s_1s_2z)^2 - k^2} + 4\pi  s_1s_2z)}\right)
e\Bigl(\epsilon_2(s_1 + \epsilon_1\epsilon_2 s_2)^2z\Bigr)
\frac{dz}{z}
\end{split}
\end{displaymath}
with $\epsilon_1, \epsilon_2 \in \{ \pm 1\}$. Repeated integration by parts (integrating only the last exponential and differentiating the rest) in the range $z \gg k^{1+\delta}$ shows that this integral is negligible unless  
$$(s_1\pm s_2)^2 \ll k^{-2\delta+ \varepsilon},$$
which in view of the above remarks implies $s_1 = s_2$. 


Next we  treat the complementary range of small $z$ by inserting  a factor $1 - V(z)$ into \eqref{mult}.  
We notice that we can compute the $\theta$-integral explicitly as a Gaussian error integral:
$$\int_0^{\pi/2} e\left(\pm \alpha \sin^2\theta\right) \sin\theta \, d\theta 
= \frac{e(\pm \alpha)}{ (2\pi  \alpha)^{1/2}}\int_0^{(2 \pi \alpha)^{1/2}} e^{\mp i t^2} dt =  \frac{1 \mp i}{4} \frac{e(\pm \alpha)}{\alpha^{1/2}}\left( 1 + O(\alpha^{-1/2})\right),$$
  see  \cite[8.253]{GR}.  
This leaves us with analyzing  
\begin{equation*} 
\begin{split}
   \int_{\Bbb{R}} \int_{\Bbb{R}}  \int_0^{\infty} & (1 - V(z)) v(x_1) v(x_2) \frac{W(n_1x_1/|q_1|^2)W(n_2x_2/|q_2|^2)}{ x_1x_2}   e\left(\pm \Bigl((s_1^2+s_2^2)z + \frac{x_1x_2  }{z}\Bigr)\right) \\ 
   &\left(\frac{1 \mp i}{4} + O\left(\frac{z^{1/2}}{(x_1x_2)^{1/2}}\right)\right) J_{\ell}\left(4\pi  s_1s_2z \right)\frac{dz}{\sqrt{z}} \,  
 e\left( \frac{x_1h_1 + x_2h_2}{|\det C|}\right)dx_1dx_ 2. 
 \end{split}
\end{equation*}
Using \eqref{bessel3}, we estimate the contribution of the error term trivially by
$$ \int_{\Bbb{R}} \int_{\Bbb{R}}  \int_{0}^{\infty}   v(x_1) v(x_2) (1 - V(z)) \frac{|W(n_1x_1/|q_1|^2)W(n_2x_2/|q_2|^2)|}{ (x_1x_2)^{3/2}} |J_{\ell} (4\pi  s_1s_2z )| dx_1\, dx_2\, dz \ll k^{-\frac{1}{2} + \delta + \varepsilon}.$$
The contribution of the main term 
is negligible, as can be seen by sufficiently many integration by parts with respect to $x_1$ or $x_2$.  Choosing $\delta$ sufficiently small completes the proof.


 
\section{Proof of Theorem \ref{thm1}}

Combining the approximate functional equation \eqref{approx} with $q=1$ with \eqref{inner}, we obtain
\begin{displaymath}
\begin{split}
\sum_{F \in \mathcal{B}_k^{(2)}} \omega_{F, k}  L(1/2, F)& = 2 \sum_{n,m} \frac{r(n) W(nm)}{(nm)^{1/2}} \sum_{F \in \mathcal{B}_k^{(2)}} c_k\frac{a_F(mI) a_F(I)}{\| F \|^2}\\
& = 2 \sum_{n,m} \frac{r(n) W(nm)}{(nm)^{1/2}}  \frac{1}{8} m^{k - 3/2} h_{mI}(I).
\end{split}
\end{displaymath}
By Lemma \ref{1}, $h_{mI}(I)$ splits into three terms. Since $\text{Aut}(I) = 8$, the diagonal contribution $m=1$ is by Mellin inversion
$$2 \sum_{n,m} \frac{r(n) W(n)}{n^{1/2}}  = \frac{2}{2\pi i} \int_{(2)} \zeta_{\Bbb{Q}(i)}(s + 1) \frac{L_{\infty}(s + 1/2)}{L_{\infty}(1/2)} (1-s^2) \frac{ds}{s}.$$
We shift the contour to $\Re s = -2 + \varepsilon$ (notice that there is no pole at $s= -1$). The remaining integral is $O(k^{-2+\varepsilon})$ by Stirling's formula, and the residue of the double pole at $s=0$ equals
$$2 L(1, \chi_{-4}) \left( \frac{\Gamma'}{\Gamma}(k-1) - \log(4\pi^2)
\right) + 2 L'(1, \chi_{-4}),$$
which gives the desired main term, since
\begin{equation}\label{digam}
\frac{\Gamma'}{\Gamma}(k-1) = \log k    + O(1/k). 
\end{equation}

Next we treat the rank 1 contribution with $T = I$, $Q = mI$. Here we must have $m \mid s$, and for given $s$ there are at most $O(s^{\varepsilon})$ choices for $U, V$. By trivial bounds using \eqref{boundW} and \eqref{bessel2} we obtain
$$ \sum_{n,m} \frac{r(n) |W(nm)|}{(nm)^{1/2}} \sum_{s, c}  \frac{c^{1/2}}{(ms)^{1/2-\varepsilon}} \Bigl|J_{\ell} \left(\frac{4\pi}{cs}\right) \Bigr| \ll k^{-100}.$$

Finally the rank 2 contribution is at most  
$$ \sum_{n,m} \frac{r(n) |W(nm)|}{(nm)^{1/2}}\sum_{\det C \not = 0} |\mathcal{J}_{\ell} (m C^{-1} C^{-\top})|.$$
By Lemma \ref{2}, the smallest eigenvalue of $m C^{-1} C^{-\top}$ is $\ll  m\| C \|^{-2}$, so that  $$ \mathcal{J}_{\ell} (m C^{-1} C^{-\top}) \ll \max_{s \ll m^{1/2}/\| C \|} |J_{\ell}(s)|$$
by trivial estimates. This gives the bound
$$  \sum_{n,m} \frac{r(n) |W(nm)|}{(nm)^{1/2}}\sum_{\det C \not = 0}   \left(\frac{m^{1/2}}{\| C \| \ell}\right)^{\ell}  \ll k^{-100}$$
by \eqref{boundW} and \eqref{bessel2} (the first of which effectively truncates $nm \ll k^{1+\varepsilon}$ up to a negligible error), and completes the proof of Theorem \ref{thm1}.\\

We conclude this section with a proof of \eqref{av}, which is a simplified version of the above computation. We have
$$\sum_{F \in \mathcal{B}_k^{(2)}} \omega_{F, k}   =  \sum_{F \in \mathcal{B}_k^{(2)}} c_k\frac{a_F(I) a_F(I)}{\| F \|^2} = \frac{1}{8}h_I(I),$$
and by \eqref{bessel2} the rank 1 and 2 contribution is trivially $O(e^{-k})$. 

\section{Symplectic exponential sums}\label{sec5}

Let 
\begin{equation}\label{defGO}
{\rm GO}_2(\Bbb{Z}) = \left\{\left(\begin{matrix} x & y\\ \mp y & \pm x \end{matrix}\right) \mid (x, y) \in \Bbb{Z}^2 \setminus \{(0, 0)\}\right\}.
\end{equation}
We denote by $\phi : \Bbb{Z}[i] \setminus \{0\} \rightarrow \Bbb{N}$ Euler's totient function on $\Bbb{Z}[i]$.  We recall the definition of $X(C)$, which is the set of (representatives of) matrices in the summation \eqref{klooster} of the Kloosterman sum.  
The aim of this section is to prove the following result.  
\begin{lemma}\label{lem4} Let   $C =  \left(\begin{smallmatrix} x & y\\ \mp y & \pm x\end{smallmatrix}\right) \in {\rm GO}_2(\Bbb{Z})$, and let $q_1, q_2$ be two coprime fundamental discriminants (possibly 1). Then
$$\mathcal{K}(C; q_1, q_2) :=   \sum_{\left(\begin{smallmatrix} A & \ast \\ C & D\end{smallmatrix}\right) \in X(C)} \sum_{\substack{\mu_1 \, (\text{{\rm mod }} [q_1, \det C])\\\mu_2 \, (\text{{\rm mod }} [q_2, \det C])}}   \chi_{q_1}(\mu_1) \chi_{q_2}(\mu_2)   e\left(\frac{\mu_1 \text{{\rm tr}}(AC^{\top}) + \mu_2 \text{{\rm tr}}(C^{\top}D)}{|\det C|}\right)$$
vanishes unless $q_1 = q_2 = 1$, in which case it equals $|\det C|^2 \phi(x+ iy)$. 
\end{lemma}

\textbf{Proof.}  The $\mu_1, \mu_2$ sum  is a product of two Gau{\ss} sums. Let us denote by $\epsilon_{q} \in \{1, i\}$ the sign of the Gau{\ss} sum associated with the character $\chi_q$. Then $\mathcal{K}(C; q_1, q_2)$    vanishes unless   $[q_1, q_2] \mid \det C$, in which case it equals
 $$\epsilon_{q_1} \epsilon_{q_2} \frac{|\det C|^2}{\sqrt{|q_1q_2|}} \sum \chi_{q_1}\left(\frac{\text{tr}(AC^{\top})}{|\det C|/|q_1|}\right) \chi_{q_2}\left(\frac{\text{tr}(C^{\top}D)}{|\det C|/|q_2|}\right),$$
 where the sum is over all
 $$\left(\begin{matrix} A & \ast \\ C & D\end{matrix}\right) \in X(C), \quad  \frac{\det C}{q_1} \mid \text{tr}(AC^{\top})  \quad  \frac{\det C}{q_2} \mid \text{tr}(C^{\top}D).   $$
We study now more carefully this set of matrices. 
  An integral matrix $M = \left(\begin{smallmatrix} A & B\\ C & D\end{smallmatrix}\right)$ is in ${\rm Sp}_4(\Bbb{Z})$ if and only if $A^{\top} C$ and $B^{\top} D$ are symmetric and $A^{\top}D - C^{\top}B = I$. This is equivalent to $M^{\top} \in {\rm Sp}_4(\Bbb{Z})$, which implies that $CD^{\top}$ is symmetric, and for $C \in {\rm GO}_2(\Bbb{Z})$ this implies that also $D^{\top}C$ is symmetric. On the other hand, 
multiplying $A^{\top}D - C^{\top}B = I$ from the left by $D^{\top} C$ gives for $C \in {\rm GO}_2(\Bbb{Z})$ that 
$$D^{\top} CA^{\top}C C^{\top}|\det C|^{-1} D - D^{\top} |\det C| B = D^{\top} C.$$
Hence if $A^{\top} C$ and $C D^{\top}$ are symmetric (and so $D^{\top} C$), then 
$D^{\top} B$ is automatically symmetric. We conclude that for $C \in {\rm GO}_2(\Bbb{Z})$ the matrix $\left(\begin{smallmatrix} A & \ast\\ C & D\end{smallmatrix}\right)$ can be extended to a symplectic matrix if and only if 
$$A^{\top} C, C D^{\top} \text{ symmetric}, \quad C(A^{\top} D - I) \equiv 0 \, (\text{mod } |\det C|).$$
Let $$C = \left(\begin{matrix} ad & bd\\ - bd &   ad\end{matrix}\right) \in {\rm GO}_2(\Bbb{Z}), \quad (a, b) = 1, d \in \Bbb{N},$$
and let $\mathcal{S}$ be the set of 2-by-2 integral symmetric matrices. 

First we determine a system of representatives of matrices $A = \left(\begin{smallmatrix} a_1 & a_2\\ a_3 &   a_4\end{smallmatrix}\right)$ modulo $\mathcal{S} \cdot C$ such that $A^{\top} C$ is symmetric. The matrix $A^{\top} C$ is symmetric if and only if 
$$a(a_3 - a_2) + b(a_1 + a_4) = 0.$$
Since $(a, b) = 1$, this is equivalent to $a_3 = a_2 - c_1 b$, $a_4 = -a_1 + c_1 a$ for some $c_1 \in \Bbb{Z}$.  
Shifting modulo $\mathcal{S} \cdot C$ with matrices $ \left(\begin{smallmatrix} 0 & 0\\ 0 &   x_4\end{smallmatrix}\right) C$, we can restrict  $c_1$ (mod $d$), and then using matrices of the form 
$$ \left(\begin{matrix} x_1 & x_2\\ x_2 &   -x_1\end{matrix}\right) C = d\left(\begin{matrix} a x_1 - b x_2 & b x_1 + a x_2\\ b x_1 + a x_2 & - a x_1 + b x_2\end{matrix}\right),$$ we can restrict $a_2$ modulo $d$ (since $(a, b) = 1$), and then $a_1$ modulo $(\det C)/d = (a^2 + b^2)d$. This is not canonical; changing the representative for $a_2$ changes simultaneously the representative for $a_1$, so we fix $a_2 \in [1, d]$ and $a_1 \in [1, (a^2 + b^2)d]$. Having used up all degrees of freedom, we conclude that a system of representatives of matrices $A$ is given by
\begin{equation}\label{A}
A = \left(\begin{matrix} a_1 & a_2\\ a_2 - c_1 b & -a_1 + c_1 a\end{matrix}\right) , \quad c_1, a_2 \in [1, d], \, a_1 \in [1, (a^2 + b^2)d].
\end{equation}
For such matrices we have $$\text{tr}(AC^{\top}) = d(a^2+b^2) c_1 = \frac{c_1 \det C}{d}.$$ 
This is divisible by $q_1^{-1}\det C$ if and only if $c_1$ is divisible by $d/(d, q_1)$. In this case 
$$\chi_{q_1}\left(\frac{\text{tr}(AC^{\top})}{|\det C|/|q_1|}\right) = \chi_{q_1}\left(\frac{c_1 |q_1|}{d}\right)$$
vanishes unless $q_1 \mid d$. 

Similarly we see that a system of   representatives of matrices $D$ modulo $C \cdot \mathcal{S}$ such that $C D^{\top}$ is symmetric  is given by
\begin{equation}\label{D}
D = \left(\begin{matrix} d_1 & d_2 +c_2b \\ d_2 & -d_1 + c_2 a\end{matrix}\right), \quad c_2, d_2 \in [1, d], \, d_1 \in [1, (a^2 + b^2)d],
\end{equation}
and we have  $q_2^{-1} \det C \mid \text{tr}(C^{\top}D)$ if and only if $c_2$ is divisible by $d/(d, q_2)$, in which case 
$$\chi_{q_2}\left(\frac{\text{tr}(C^{\top}D)}{|\det C|/|q_2|}\right) = \chi_{q_2}\left(\frac{c_2 |q_2|}{d}\right)$$
vanishes unless $q_2 \mid d$. We summarize
$$\mathcal{K}(C; q_1, q_2) = \delta_{[q_1, q_2] \mid d} \epsilon_{q_1} \epsilon_{q_2} \frac{|\det C|^2}{\sqrt{|q_1q_2|}} \sum  \chi_{q_1}\left(\frac{c_1 |q_1|}{d}\right) \chi_{q_2}\left(\frac{c_2 |q_2|}{d}\right)  $$
where the sum is over all pairs $(A, D)$ as in \eqref{A} and \eqref{D} satisfying 
\begin{equation}\label{divi}
\frac{d}{q_1} \mid c_1, \quad \frac{d}{q_2} \mid c_2, \quad C(A^{\top} D - I) \equiv 0 \, (\text{mod } \det C).
\end{equation}
 The latter condition implies in particular $A^{\top} D \equiv I$ (mod $d$), so that $\det(A^{\top} D) \equiv 1$ (mod $d$). In particular, $A$ and $D$ are invertible modulo $d$, and we obtain $D \equiv A^{-\top} = (A^{\top})^{\text{adj}}\, \overline{\det A}$ (mod $d$). Considering the trace and the difference of the off-diagonal entries of this congruence, we conclude $c_2 a \equiv c_1 a \, \overline{\det A} \, (\text{mod } d)$ and $c_2b \equiv  c_1b \, \overline{\det A} \, (\text{mod } d)$. Since $(a, b) = 1$, this implies $c_2 \equiv c_1 \, \overline{\det A} \, (\text{mod } d)$, and in particular $(c_2, d) = (c_1, d)$. Hence $\mathcal{K}(C; q_1, q_2)$ contains a subsum
$$\sum_{\substack{A \, (\text{mod d})\\ c_1 \equiv 0 \, (\text{mod } d/(q_1, q_2))}} \chi_{q_1}\left(\frac{c_1   |q_1|}{d}\right) \chi_{q_2}\left(\frac{c_1 |q_2| \det A}{d}\right),$$
 which obviously vanishes unless\footnote{In fact, one can show that it vanishes unless $q_1 = q_2 \in \{1, -4\}$.} $|q_1| = |q_2|$. Since $(q_1, q_2) = 1$, we are left with analyzing the case $q_1 = q_2 = 1$, where we can assume $c_1 = c_2 = 0$. Putting $\gamma = (a^2 + b^2) d$, this leaves us with counting matrices 
$$A = \left(\begin{matrix} a_1 & a_2\\ a_2   & -a_1 \end{matrix}\right), \quad D = \left(\begin{matrix} d_1 & d_2   \\ d_2 & -d_1  \end{matrix}\right)  , \quad  a_2, d_2 \in [1, d], \, a_1, d_1 \in [1, \gamma]$$
satisfying $\left(\begin{smallmatrix} a & b\\ -b & a \end{smallmatrix}\right) (A^{\top} D - I) \equiv 0 \, (\text{mod } \gamma).$ 
  In order to avoid problems with well-definedness we write this as 
  $$\mathcal{K}(C;1, 1) =   \frac{|\det C|^2}{(a^2+b^2)^2} \#\Bigl\{a_1, a_2, d_1, d_2 \, (\text{mod } \gamma)  \mid \left(\begin{matrix} a & b\\ -b & a \end{matrix}\right) (A^{\top} D - I) \equiv 0 \, (\text{mod } \gamma) \Bigr\},$$
where
$ A = \left(\begin{smallmatrix} a_1 & a_2\\ a_2 & -a_1\end{smallmatrix}\right)$, $D = \left(\begin{smallmatrix} d_1 & d_2\\ d_2 & -d_1\end{smallmatrix}\right).$ 
This can be conveniently rephrased as
\begin{displaymath}
\begin{split}
\mathcal{K}(C;1, 1)  & = \frac{|\det C|^2}{(a^2+b^2)^2}\#\bigl\{\alpha, \delta \in \Bbb{Z}[i]/c\Bbb{Z}[i] \mid \alpha \delta \equiv 1 \, (\text{mod } ad + bdi)\bigr\} = |\det C|^2 \phi(ad + bd i).
\end{split}
\end{displaymath}
The same analysis works if $C = \left(\begin{smallmatrix} a d & b d\\ bd & - a d\end{smallmatrix}\right)$ has negative discriminant, and the proof is complete.

\section{Proof of Theorem \ref{thm2}}\label{sec6}

In this section we regard   $q_1, q_2$ as fixed, but it is clear that all implied constants depend polynomially on these quantities. 
We start similarly as in the proof of Theorem \ref{thm1} with the approximate functional equation getting
\begin{displaymath}
\begin{split}
& \sum_{F \in \mathcal{B}_k^{(2)}} \omega_{F, k} L(1/2, F \times \chi_{q_1})L(1/2, F \times \chi_{q_2})\\
& = 4 \sum_{n_1,m_1, n_2, m_3} \frac{r_{q_1} (n_1)r_{q_2} (n_2) W(n_1m_1/|q_1|^2)W(n_2m_2/|q_2|^2) \chi_{q_1}(m_1) \chi_{q_2}(m_2) }{(n_1m_1n_2m_2)^{1/2}}\\
& \quad\quad\quad\quad\quad\quad \times   \sum_{F \in \mathcal{B}_k^{(2)}} c_k\frac{a_F(m_1I) a_F(m_2I)}{  \| F \|^2}\\
& = 4\sum_{n_1,m_1, n_2, m_2} \frac{r_{q_1} (n_1)r_{q_2} (n_2) W(n_1m_1/|q_1|^2)W(n_2m_2/|q_2|^2) \chi_{q_1}(m_1) \chi_{q_2}(m_2)}{(n_1m_1n_2m_2)^{1/2}}    \frac{1}{8} m_2^{k - 3/2} h_{m_2I}(m_1I) .
\end{split}
\end{displaymath}
For notational simplicity let us write $r_1$ and $r_2$ for $r_{q_1}$ and $r_{q_2}$. The diagonal term equals
$$ 4\sum_{n_1,n_2, m} \frac{r_{1} (n_1)r_{2} (n_2) W(n_1m_1/|q_1|^2)W(n_2m_2/|q_2|^2) \chi_{q_1q_2}(m)}{(n_1n_2m^2)^{1/2}}.  $$
By 
 double Mellin inversion this is
\begin{equation}\label{integrand}
\begin{split}
&= \frac{4}{(2\pi i)^2} \int_{(2)} \int_{(2)} L(s+1, \chi_{q_1}) L(s+1, \chi_{-4q_1 })L(t+1, \chi_{q_2}) L(t+1, \chi_{-4q_2 })   L(s+t+1, \chi_{q_1 q_2}) \\
& \quad\quad\quad\quad\quad \times \frac{L_{\infty}(s+1/2)}{L_{\infty}(1/2)}  \frac{L_{\infty}(t+1/2)}{L_{\infty}(1/2)}(1-s^2)(1-t^2) |q_1|^{2s} |q_2|^{2t} \frac{ds\, dt}{st}.
\end{split}
\end{equation}
We shift the $s$-contour to $\Re s = -1$, picking up a  pole at $s=0$ of order 1 or 2. For the latter, we shift the $t$-contour to $\Re t = -1$, picking up a   pole at $t=0$ of order at most 4; 
the remaining integral is $O(k^{-1+\varepsilon})$. 


For the former, we shift the $t$-contour to $\Re t = - 1 $, picking up a possible pole at $t = -s$ (since $(q_1, q_2) = 1$, this can happen only if $q_1 = q_2= 1$) and a   pole at $t=0$. The latter as well as the remaining integral contributes $O(k^{-1+\varepsilon})$; 
if $q_1 = q_2 =1$, the former equals
\begin{equation}\label{polar}
-\frac{4}{2\pi i}\int_{(-1+\varepsilon)} \Gamma(s+1)\zeta_{\Bbb{Q}(i)}(s+1) \Gamma(1-s)\zeta_{\Bbb{Q}(i)}(1-s)  \frac{\Gamma(k-1+s)\Gamma(k-1-s)}{\Gamma(k-1)^2}(1-s^2)^2 \frac{ds}{s^2}.
\end{equation}
We can evaluate this term as follows: we move the contour to $\Re s = \varepsilon$  and truncate it at $|\Im s| \leq k^{\varepsilon}$ at the cost of a negligible error.  
By Stirling's formula one obtains
\begin{equation*}
\frac{\Gamma(k-1+s)\Gamma(k-1-s)}{\Gamma(k-1)^2} = 1 + O\left( \frac{(\Im s)^2}{ k}\right)
\end{equation*}
for $\Im s \ll k^{1/3}$, say. Adding back the truncated contour shows that the integral in questions equals
$$-\frac{4}{2\pi i}\int_{(\varepsilon)} \Gamma(s+1)\zeta_{\Bbb{Q}(i)}(s+1) \Gamma(1-s)\zeta_{\Bbb{Q}(i)}(1-s)  (1-s^2)^2 \frac{ds}{s^2} + O(k^{-1+\varepsilon}),$$
and   the main term is a constant independent of $k$. However, this maneuver is not necessary, since the term \eqref{polar} will be cancelled by another term in a moment. \\

Next we turn to the rank 1 contribution. Here we must have $[m_1, m_2] \mid s$, and by trivial estimates we obtain
\begin{displaymath}
\begin{split}
&  \sum_{n_1,m_1, n_2, m_2} \frac{  |W(n_1m_1/|q_1|^2)W(n_2m_2/|q_2|^2)|}{n_1n_2 ( m_1 m_2)^{1/2}}  \sum_{s, c} \frac{c^{1/2}}{([m_1, m_2] s)^{1/2 - \varepsilon}} \Bigl|J_{\ell} \left(\frac{4\pi m_1m_2   }{[m_1, m_2] c s}\right)\Bigr|. 
\end{split}
\end{displaymath}
Here we can truncate the $m_1, m_2, n_1, n_2$-sum at 
$n_1m_1, n_2m_2 \ll k^{1+\varepsilon}$ 
at the cost of a negligible error, and then by \eqref{bessel2} the $c, s$-sum at $cs \leq  50(m_1,m_2)/k$, again at the cost of a negligible error. 
By trivial estimates we obtain the bound
$$  k^{\varepsilon} \sum_{\substack{n_1dm_1 \leq   k^{1+\varepsilon}\\ n_2dm_2 \leq   k^{1+\varepsilon} \\ sc \ll d /k} }  \frac{ c^{1/2}}{n_1n_2m_1m_2 s^{1/2} d^{3/2}}  \ll  k^{-1/2+\varepsilon}, 
$$
which majorizes the above error terms in the residue computation. \\

It remains to treat the rank 2 contribution
\begin{displaymath}
\begin{split}
4 \pi^2 \sum_{n_1,m_1, n_2, m_2} &\frac{r_1(n_1)r_2(n_2) W(n_1m_1/|q_1|^2)W(n_2m_2/|q_2|^2) \chi_{q_1}(m_1)\chi_{q_2}(m_2)}{(n_1m_1n_2m_2)^{1/2}} \\
& \times \sum_{\det C \not = 0} \frac{K(m_2I, m_1I; C)}{|\det C|^{3/2}} \mathcal{J}_{\ell}(m_1m_2 C^{-1} C^{-\top}). 
\end{split}
\end{displaymath}
By the decay \eqref{boundW} of the weight function we can again truncate the $m_1, m_2, n_1, n_2$-sum at $n_1m_1, n_2m_2 \ll k^{1+\varepsilon}$ at the cost of a negligible error. Then by the same argument as in the proof of Theorem \ref{thm1} we can truncate the $C$-sum at $\| C \| \ll k^{\varepsilon}$  at the cost of a negligible error by the rapid decay \eqref{bessel2} of the Bessel function for arguments less than the index.  Having truncated the $C$-sum, we can complete the $m_1, m_2, n_1, n_2$-sum at the cost of a negligible error by \eqref{boundW} and apply Poisson summation in $m_1, m_2$ split into residue classes modulo $[q_1, \det C]$ and $[q_2, \det C]$ respectively.  This gives
\begin{displaymath}
\begin{split}
4\pi^2 \sum_{n_1,  n_2 } &\frac{r_1(n_1)r_2(n_2) }{(n_1 n_2 )^{1/2}}  \sum_{\substack{\det C \not = 0\\ \| C \| \ll k^{\varepsilon}}}   \sum_{\substack{\mu_1  \, (\text{mod } [q_1, \det C])\\ \mu_2 \, (\text{mod } [q_2, \det C])}} \chi_{q_1}(\mu_1) \chi_{q_2}(\mu_2) \\
& \times \sum_{h_1, h_2 \in \Bbb{Z}} \frac{K(\mu_2I, \mu_1I; C)}{[q_1, \det C] [q_2, \det C] |\det C|^{3/2}}e\left(- \frac{\mu_1h_1 + \mu_2h_2}{|\det C|}\right) \Psi_{n_1, n_2, q_1, q_2}(C; h_1, h_2), \\
\end{split}
\end{displaymath}
 where $\Psi_{n_1, n_2, q_1, q_2}(C; h_1, h_2) = \Psi(C; h_1, h_2)$ was defined in \eqref{defPsi}. 
 Integrating by parts sufficiently often with the help of \eqref{diff}, we can truncate the $h_1, h_2$-sum at $k^{\varepsilon}$ at the cost of a negligible error. By Lemma \ref{psi}, the previous display equals
 \begin{equation}\label{equals}
 \begin{split}
4\pi^2 \sum_{n_1,  n_2 }&  \frac{r_1(n_1)r_2(n_2) }{(n_1 n_2 )^{1/2}}  \sum_{\substack{C \in {\rm GO}_2(\Bbb{Z})\\ \| C \| \ll k^{\varepsilon}}}     \frac{\Psi(C; 0, 0)  }{[q_1, \det C] [q_2, \det C] |\det C|^{3/2}}  \\
&\times  \sum_{\substack{\mu_1  \, (\text{mod } [q_1, \det C])\\ \mu_2 \, (\text{mod } [q_2, \det C])}} \chi_{q_1}(\mu_1) \chi_{q_2}(\mu_2)  K(\mu_2I, \mu_1I; C) + O(k^{-1/2 + \varepsilon}),
 \end{split}
 \end{equation}
where ${\rm GO}_2(\Bbb{Z})$ was defined in \eqref{defGO} and 
equals precisely the set of invertible integral 2-by-2 matrices such that  $C^{\top} C$ has two identical eigenvalues (i.e.\ is a multiple of the identity).  
For $C \in {\rm GO}_2(\Bbb{Z})$ we have $C^{-1} = C^{\top}|\det C|^{-1}$, so that by the definition \eqref{klooster} of the Kloosterman sum the double sum over $\mu_1, \mu_2$ equals 
 $\mathcal{K}(C, q_1, q_2)$ as defined in Lemma \ref{lem4}. 
Hence the main term in \eqref{equals} equals
\begin{equation}\label{equals1}
\begin{split}
8\pi^2 \delta_{q_1= q_2 = 1} \sum_{n_1,  n_2 } \frac{r_1(n_1)r_2(n_2) }{(n_1 n_2 )^{1/2}}  \sum_{\substack{\gamma \in \Bbb{Z}[i] \setminus \{0\} \\ |\gamma|^2  \ll k^{\varepsilon}}} \frac{\phi(\gamma)}{|\gamma|^3}  \int_{\Bbb{R}} \int_{\Bbb{R}} \frac{W(n_1x_1) W(n_2x_2)}{\sqrt{x_1x_2}} \mathcal{J}_{\ell}\left(\frac{x_1x_2}{|\gamma|^2} I  \right) 
dx_1dx_2,
 \end{split}
 \end{equation}
and by \eqref{bessel2} and \eqref{boundW} we can complete the $\gamma$-sum at the cost of a negligible error. By \eqref{mellin} and Mellin inversion we have
\begin{displaymath}
\begin{split}
&\sum_{ \gamma \in \Bbb{Z}[i] \setminus \{0\}} \frac{\phi(\gamma)}{|\gamma|^3}  \mathcal{J}_{\ell}\left(\frac{x_1x_2}{|\gamma|^2} I  \right)\\
&  = \int_0^{\pi/2} \int_{(2)}  \sum_{ \gamma \in \Bbb{Z}[i] \setminus \{0\}} \frac{\phi(\gamma)}{|\gamma|^{3+s}}   \frac{(4\pi \sqrt{x_1x_2} \sin\theta)^s \Gamma(k - \frac{3+s}{2}) \Gamma(\frac{1 + s}{2})}{2\pi^{1/2}  \Gamma(k + \frac{s-1}{2}) \Gamma(1+\frac{s}{2})}\frac{ds}{2\pi i}  \sin \theta \, d\theta.
\end{split}
\end{displaymath}
Notice that the right hand side  is absolutely convergent. 
We compute the $\theta$-integral explicitly \cite[3.621.1]{GR}
$$ \int_0^{\pi/2} (\sin\theta)^{s+1} d\theta = \frac{\pi^{1/2} \Gamma(1 + s/2)}{2 \Gamma((3+s)/2)},$$
 as well as the $\gamma$-sum (noting that $\Bbb{Z}[i]$ has 4 units), so that 
$$\sum_{ \gamma \in \Bbb{Z}[i] \setminus \{0\}} \frac{\phi(\gamma)}{|\gamma|^3}  \mathcal{J}_{\ell}\left(\frac{x_1x_2}{|\gamma|^2} I  \right) =  2 \int_{(2)}  \frac{\zeta_{\Bbb{Q}(i)}((s+1)/2)}{\zeta_{\Bbb{Q}(i)}((s+3)/2)} \frac{(4\pi \sqrt{x_1x_2} )^s \Gamma(k - \frac{3+s}{2})}{ (1+s)  \Gamma(k + \frac{s-1}{2}) }\frac{ds}{2\pi i}.$$
Computing now the $x_1, x_2$-integral,   we can re-write \eqref{equals1} (up to a negligible error) for $q_1 = q_2 = 1$ as
\begin{displaymath}
\begin{split}
& 16\pi^2 \int_{(2)}  \sum_{n_1,  n_2 } \frac{r(n_1)r(n_2) }{(n_1 n_2 )^{1 + s/2}}  \frac{\zeta_{\Bbb{Q}(i)}((s+1)/2)}{\zeta_{\Bbb{Q}(i)}((s+3)/2)} \Bigl(\frac{ L_{\infty}(1 + s/2) (1 - (\frac{s+1}{2})^2)}{L_{\infty}(1/2) (s+1)/2}\Bigr)^2 \frac{(4\pi  )^s \Gamma(k - \frac{3+s}{2})}{ (1+s)  \Gamma(k + \frac{s-1}{2}) }\frac{ds}{2\pi i}\\
=  &8\pi^2 \int_{(2)}  \zeta_{\Bbb{Q}(i)}\Bigl(\frac{s+3}{2} \Bigr)  \zeta_{\Bbb{Q}(i)}\Bigl(\frac{s+1}{2} \Bigr)    \frac{\Gamma(k - \frac{s+3}{2})\Gamma(k + \frac{s-1}{2}) \Gamma(\frac{s+1}{2}) \Gamma(\frac{s+3}{2})}{\pi^{2+s}(1+s)^2 \Gamma(k-1)^2} \left(1 - \Bigl(\frac{s+1}{2}\Bigr)^2\right)^2 \frac{ds}{2\pi i}.
 \end{split}
 \end{displaymath}
 Changing variables and applying the functional equation of the Dedekind zeta-function, this equals
 \begin{displaymath}
 \begin{split}
 & 4\pi^2 \int_{(3/2)}  \zeta_{\Bbb{Q}(i)}(s+1)  \zeta_{\Bbb{Q}(i)}(s)\frac{\Gamma(k -1-s)\Gamma(k -1+s) \Gamma(s) \Gamma(s+1)}{\pi^{1+2s}s^2 \Gamma(k-1)^2} (1 - s^2)^2\frac{ds}{2\pi i}\\
= & 4 \int_{(3/2)}  \zeta_{\Bbb{Q}(i)}(s+1)  \zeta_{\Bbb{Q}(i)}(1-s)\frac{ \Gamma(k -1-s)\Gamma(k -1+s) \Gamma(1-s) \Gamma(s+1)}{ s^2 \Gamma(k-1)^2} (1 - s^2)^2 \frac{ds}{2\pi i}.
 \end{split}
 \end{displaymath}
 This cancels precisely the term \eqref{polar} (notice that the residue at $s=0$ vanishes) and 
completes the proof of Theorem \ref{thm2}. 

\section{The level aspect}

In this section we sketch a proof of Theorem \ref{thm3}. Let $N \equiv 3$ (mod 4) be a (large) prime. The group $\Gamma = \Gamma_0^{(2)}(N)$ has index $\asymp N^{3}$ in ${\rm Sp}_4(\Bbb{Z})$. In particular, by our normalization \eqref{inner1} of the inner product, an $L^2$-normalized form $F$ of level 1 viewed as an oldform of level $N$ has norm $\asymp N^3$. 

 If  $F$ is a Siegel modular \emph{new}form for $\Gamma_0^{(2)}(N)$ that is not a lift, then  $L(s, F)$ has conductor $N^2$ in all cases where $\omega_{F, k} \not= 0$. Indeed\footnote{this argument   was communicated to the author by A.\ Saha},  the local representation attached to   $F$ at the prime  $N$ is one of  IIa, IIIa, Vb/c, VIa, VIb, see \cite[p.\ 13]{DPSS}.  But IIa, Vb/c, and VIa have no corresponding local Bessel model, so $a_F(I) = 0$, and the remaining types have conductor $N^2$.  For such $F$, we obtain the approximate functional equation 
$$\omega_{F, k} L(1/2, F) = 2 c_k \sum_{n,m} \frac{r(n) W(nm/N)}{(nm)^{1/2}}   \frac{a_F(mI) a_F(I)}{\| F \|^2}$$
with $W$ as in \eqref{approx}.
Both the right hand side and the left hand side make sense also for oldforms and lifts (but they do not have to be identical in this case). Let $\mathcal{B}_k^{(2)}(N)$ be a Hecke basis for $\Gamma$. For the $O(1)$ oldforms in $\mathcal{B}_k^{(2)}(N)$ of level 1, the left-hand side is $O(N^{-3})$ and the right hand side is 
$$\ll N^{-3+\varepsilon} \frac{|a_F(I)|^2}{\| F \|^2_{{\rm Sp}_4(\Bbb{Z})}} \sum_{m \ll N^{1+\varepsilon}} \frac{|\lambda_F(m)|}{m^{1/2}} \ll N^{-2+ \varepsilon}$$
by the trivial bound $\lambda_F(m) \ll m^{1/2+\varepsilon}$ (which is sharp for lifts). For each of the $O(N)$ lifts in $\mathcal{B}_k^{(2)}(N)$, the left hand side is $N^{-3 +\frac{1}{4} + \varepsilon}$ by the convexity bound for ${\rm GL}_2$ $L$-functions and the fact that B\"ocherer's conjecture is known in this case so that $\omega_{F, k} \ll N^{-3+\varepsilon}$. The right hand side, by the same argument as before, is $O(N^{-2+\varepsilon})$. Hence we conclude that
\begin{displaymath}
\begin{split}
\sum_{F \in \mathcal{B}_k^{(2)}(N)} \omega_{F, k}  L(1/2, F)& = 2 \sum_{n,m} \frac{r(n) W(nm/N)}{(nm)^{1/2}} \sum_{F \in \mathcal{B}_k^{(2)}(N)} c_k\frac{a_F(mI) a_F(I)}{\| F \|^2} + O(N^{-1+\varepsilon}).
\end{split}
\end{displaymath}
We can now proceed as in the proof of Theorem \ref{thm1}, using the extra divisibility conditions as in Remark \ref{rem2}. In the following $k \geq 6 $ is regarded as fixed. The diagonal contribution is 
\begin{displaymath}
\frac{2}{2\pi i} \int_{(2)} \zeta_{\Bbb{Q}(i)}(s+1) \frac{L_{\infty}(s+1/2)}{L_{\infty}(1/2)} (1-s^2) N^s \frac{ds}{s} = 2 L(1, \chi_{-4})\log N + C(k) + O(N^{-2+\varepsilon}),
\end{displaymath}
where
$$C(k) = 2L(1, \chi_{-4}) \left(\frac{\Gamma'}{\Gamma}(k-1) + \log\frac{1}{4\pi^2}\right) + 2L'(1, \chi_{-4}).$$
The rank 1 contribution is trivially bounded by 
$$ \sum_{n,m} \frac{r(n) W(nm/N)}{(nm)^{1/2}} \sum_{s, c}  \frac{(Nc)^{1/2}}{(ms)^{1/2-\varepsilon}} \Bigl|J_{\ell} \left(\frac{4\pi}{Ncs}\right) \Bigr| \ll N^{\varepsilon} \sum_{N \mid c} \frac{1}{c^{\ell - 1/2}} \ll N^{2 - k + \varepsilon}. $$
The rank 2 contribution is trivially bounded by
$$ \sum_{n,m} \frac{r(n) W(nm)}{(nm)^{1/2}}\sum_{\substack{\det C \not = 0\\ N \mid C}} |\mathcal{J}_{\ell} (m C^{-1} C^{-\top})| \ll N^{-\ell + \varepsilon} \sum_{m \ll N^{1+\varepsilon}}  m^{\frac{\ell}{2} - \frac{1}{2}} \ll N^{\frac{7}{4} - \frac{k}{2}+\varepsilon}.$$

\end{document}